\documentstyle[sprocl]{article}

\def\Journal#1#2#3#4{{#1} {\bf #2}, #3 (#4)}

\def\be{\begin{equation}}
\def\ee{\end{equation}}
\def\bea{\begin{eqnarray}}
\def\eea{\end{eqnarray}}

\begin{document}

\hfill q-alg/9710009

\title{QUANTUM ORTHOGONAL CALEY-KLEIN GROUPS AND ALGEBRAS}

\author{N.A. GROMOV, I.V. KOSTYAKOV, V.V. KURATOV}

\address{Department of Mathematics, Komi Science Center, Ural
 Division,\\
 Russian Academy of Science, Syktyvkar, 167610, Russia}




\maketitle\abstracts{
 The extension of FRT quantization theory
for the nonsemisimple CK groups is suggested. The quantum
orthogonal CK groups are realized as the Hopf algebras
of the noncommutative functions over an associative algebras
with nilpotent commutative generators. The quantum CK algebras
are obtained as the dual objects to the corresponding
quantum groups.}

\section{Introduction}

The quantization theory of the simple groups and algebras Lie was
developed by Faddeev-Reshetikhin-Takhtadjan (FRT) \cite{Fad--89}.
In group theory there is a remarkable set of groups, namely
the motion groups of n-dimensional spaces of constant curvature
or the orthogonal Cayley-Klein (CK) groups.
The well known
Euclidean $ E(n), $ Poincare $ P(n), $ Galileian $ G(n) $ and
other nonsemisimple groups are in the set of CK groups.
   The principal proposal for quantization of all CK groups
    is to regard them as the groups
over an  algebra $ {\bf D} $ with nilpotent
commutative generators
and the corresponding quantum CK
groups as the algebra of noncommutative functions over $ {\bf D}. $

     Algebra $ {\bf D}_{n}(\iota;{\bf C}) $ is defined as an associative
  algebra with unit and {\it nilpotent} generators
   $ \iota_{1},\ldots,\iota_{n},\; \iota_{k}^{2}=0,\; k=1,\ldots,n $
 with {\it commutative} multiplication
     $ \iota_{k}\iota_{m}=\iota_{m}\iota_{k},\; k\neq{m}. $
  The general element of $ {\bf D}_{n}(\iota;{\bf C}) $ has the form
$
a=a_0 + \sum_{p=1}^{n} \sum_{k_1<...<k_p}
a_{k_1...k_p}\iota_{k_1}...\iota_{k_p}, \quad
a_0,a_{k_1...k_p}\in{\bf C}.
$
For $ n=1 $ we have
$ {\bf D}_{1}(\iota_{1};{\bf C})\ni{a=a_{0}+a_{1}\iota_{1}}, $
i.e. dual (or Study) numbers.
For $ n=2 $ the general element of
$ {\bf D}_{2}(\iota_{1},\iota_{2};{\bf C}) $
is written as follows:
$ a=a_{0}+a_{1}\iota_{1}+a_{2}\iota_{2}+a_{12}\iota_{1}\iota_{2}. $

\section{ Orthogonal CK groups $ SO(N;j;{\bf R}) $  }

Let us regard a vector space ${\bf R}_{N}(j) $ over
$ {\bf D}_{N-1}(j;{\bf R}) $
with Cartesian coordinates
$ x(j)=(x_1,J_{12}x_2,\ldots, J_{1,N}x_{N})^t, \; x_k\in{\bf R},
 \;k=1,\ldots,N $
and quadratic form
$
x^t(j)x(j) = x_1^2 + \sum_{k=2}^{N}J_{1k}^2x_k^2,
$
where
 $
 J_{\mu\nu}=\prod_{r=\mu}^{\nu -1}j_{r}, \;
     \mu < \nu, \; J_{\mu\nu}=1,\; \mu \geq \nu, \; j_{r}=1,\iota_{r},i.
  $

 Orthogonal CK groups $ SO(N;j;{\bf R}) $
 are defined as the
 set of transformations of $ {\bf R}_{N}(j) $ leaving invariant
 $ x^t(j)x(j) $ and are
 realized in the
 Cartesian basis as the matrix groups over ${\bf D}_{N-1}(j;{\bf R}) $
 with the help of the {\it special} matrices
     \begin{equation}
   (A(j))_{kp}=\tilde{J}_{kp}a_{kp},  \quad a_{kp}\in{\bf R},
 \quad
  \tilde{J}_{kp}=J_{kp},\; k<p,   \quad    \tilde{J}_{kp}=J_{pk},\;
k\geq{p},
    \label{16}
    \end{equation}
 These matrices  act on  vectors  $ x(j)\in{\bf R}_{N}(j) $
 by matrix multiplication and
  are   satisfied   the  following $ j$-orthogonality
relations:
      $
    A(j)A^{t}(j)=A^{t}(j)A(j)=I.
      $

     One of the solutions of matrix equation
$ DC_0D^t=I, \; (C_0)_{ik}=\delta_{ik'}, \; k'=N+1-k $
provide the similarity transformation
      $
 B(j)=D^{-1}A(j)D
      $
which give the realization of
$ SO_{q}(N;{\bf C}) $  in a new ("symplectic") basis
with invariant quadratic form
      $
x^t(j)C_0x(j)
      $
and  the additional relations of $ j$-orthogonality
      $
   B(j)C_{0}B^{t}(j)=B^{t}(j)C_{0}B(j)=C_{0}.
      $

\section{Quantum groups and algebras}

     We shall regard the quantum deformations of the contracted CK
 groups, i.e. $ j_{k}=1,\iota_{k}. $
 We shall start with the
   $ {\bf D}\langle t_{ik} \rangle  $  ---  the  algebra   of   noncommutative
 polynomials of $ N^2 $ variables $ t_{ik},\; i,k=1,\ldots,N $  over
 the  algebra $ {\bf D}_{N-1}(j). $ In addition we shall transform
 the deformation parameter $ q=\exp{z} $ as follows:
   $
     z=Jv, \; J\equiv J_{1N}=\prod_{k=1}^{N-1}j_{k},
   $
 where $ v $ is the new deformation parameter.

      In "symplectic" basis the quantum CK group $ SO_{v}(N;j;{\bf C}) $
 is produced by the generating matrix
$ T(j)\in{M_{N}({\bf D} \langle t_{ik} \rangle ) } $
 equal to $ B(j) $  for $ q=1. $ The non\-com\-mu\-ta\-ti\-ve entries of
 $ T(j) $ obey the commutation relations
\begin{equation}
 R_v(j)T_1(j)T_2(j)=T_2(j)T_1(j)R_v(j).
\label{21}
\end{equation}
 and the additional relations of $ (v,j)$-orthogonality
     \begin{equation}
   T(j)C(j)T^{t}(j)=T^{t}(j)C(j)T(j)=C(j),
    \label{22}
    \end{equation}
 where lower triangular R--matrix $ R_{v}(j) $ and $ C(j) $ are
 obtained from $ R_q $ and $ C, $ respectively, by substitution
 $ Jv $ instead of $ z: $
$
     R_{v}(j)=R_{q}(z \rightarrow Jv), \quad
     C(j)=C(z \rightarrow Jv).
     $
 Then the quotient
     \begin{equation}
SO_{v}(N;j;{\bf C})= {\bf D} \langle t_{ik} \rangle \big/ (\ref{21}),(\ref{22})
    \label{24}
    \end{equation}
is Hopf algebra with the  coproduct $ \Delta, $
counit $ \epsilon $ and antipode $ S: $
\begin{equation}
\Delta T(j)=T(j) \dot {\otimes}T(j),\quad \epsilon (T(j))=I, \quad
 S(T(j))=C(j)T^t(j)C^{-1}(j).
\label{25}
\end{equation}

 By FRT quantization  theory \cite{Fad--89} the dual space
   $ Hom(SO_{v}(N;j;{\bf C}),{\bf C}) $ is an algebra with the multiplication
 induced by coproduct $ \Delta $ in $ SO_{v}(N;j;{\bf C}) $
     \begin{equation}
     l_{1}l_{2}(a)=(l_{1}\otimes l_{2})(\Delta (a)),
    \label{30}
    \end{equation}
  $ l_{1},l_{2}\in{ Hom(SO_{v}(N;j;{\bf C}),{\bf C})}, \quad a\in{SO_{v}(N;j)}. $
 Let us formally introduce $ N \times N $ upper $ (+) $ and lower
 $ (-) $ triangular matrices $ L^{(\pm)}(j) $ as follows: it is
 necessary to put $ j_{k}^{-1} $ in the nondiagonal matrix elements
  of $ L^{(\pm)}(j), $ if there is the parameter $ j_k $ in the
  corresponding matrix element of $ T(j). $ For example, if
  $ (T(j))_{12}=j_{1}t_{12}+j_{2}\tilde{t}_{12}, $ then
  $ (L^{(+)}(j))_{12}=j_{1}^{-1}l_{12}+j_{2}^{-1}\tilde{l}_{12}. $
  Formally the matrices   $   L^{(\pm)}(j)   $   are   not   defined  for
  $ j_{k}=\iota_{k}, $ since $ \iota_{k}^{-1} $ do not exist, but
  if we set an action of the matrix functionals $ L^{(\pm)}(j) $
  on the elements of $ SO_{v}(N;j;{\bf C}) $ by the duality relation
     \begin{equation}
    \langle L^{(\pm)}(j),T(j) \rangle = R^{(\pm)}(j),
    \label{31}
    \end{equation}
where
   $
     R^{(+)}(j)=PR_{v}(j)P, \;  R^{(-)}(j)= R_{v}^{-1}(j),\;
     Pu \otimes w = w \otimes u,
     $
then we shall have well defined expressions even for
$ j_{k}=\iota_k. $

     The elements of $ L^{(\pm)}(j) $ satisfy the commutation
 relations
     \begin{eqnarray}
     R^{(+)}(j)L_{1}^{(\sigma)}(j)L_{2}^{(\sigma)}(j) & = &
 L_{2}^{(\sigma)}(j)L_{1}^{(\sigma)}(j)R^{(+)}(j), \nonumber  \\
     R^{(+)}(j)L_{1}^{(+)}(j)L_{2}^{(-)}(j) & = &
L_{2}^{(-)}(j)L_{1}^{(+)}(j)R^{(+)}(j),  \quad
     \sigma = \pm
    \label{33}
    \end{eqnarray}
and additional relations
     \begin{eqnarray}
  L^{(\pm)}(j)C^{t}(j)L^{(\pm)}(j) & = & C^{t}(j),  \nonumber \\
L^{(\pm)}(j)(C^{t}(j))^{-1}L^{(\pm)}(j)  & = & (C^{t}(j))^{-1}, \nonumber \\
l_{kk}^{(+)}l_{kk}^{(-)}=l_{kk}^{(-)}l_{kk}^{(+)}=1, & &
     l_{11}^{(+)}\ldots l_{NN}^{(+)}=1, \; k=1,\ldots ,N.
    \label{34}
    \end{eqnarray}
 An algebra  $  so_{v}(N;j;{\bf C})=\{I,L^{(\pm)}(j)\} $ is called quantum CK
algebra and
is Hopf algebra with the following  coproduct $ \Delta, $
antipode $ S $ and  counit $ \epsilon : $
$
 \Delta L^{(\pm)}(j)=L^{(\pm)}(j) \dot {\otimes} L^{(\pm)}(j), \quad
S(L^{(\pm)}(j))=C^{t}(j)(L^{(\pm)}(j))^{t}(C^{t}(j))^{-1}, \\
\epsilon (L^{(\pm)}(j))=I.
$

 It is possible to show that algebra $ so_{v}(N;j;{\bf C}) $ is isomorphic
 with the quantum deformation \cite{Val--95} of the universal enveloping algebra
 of the CK algebra $ so(N;j;{\bf C}), $ which may be obtained from the
 orthogonal algebra $ so(N;{\bf C}) $ by contractions \cite{Mon--90}.


\section*{References}

\begin{thebibliography}{99}

\bibitem{Fad--89}N.Yu.  Reshetikhin, L.A. Takhtadjan and L.D.
 Faddeev, \Journal{\em Algebra and Anaysis}{1}{178}{1989} (in Russian).

\bibitem{Mon--90}N.A. Gromov,  {\em Contractions and Analytical
 Continuations of the Classical Groups: Unified Approach.}
(Komi Science Center, Syktyvkar, 1990) (in Russian).

\bibitem{Val--95}A. Ballesteros, F.G. Herranz, M.A. del Olmo and M. Santander,
\Journal{\em Lett.Math.Phys.}{33}{273}{1995}.





\end{thebibliography}
\end{document}